\documentclass[12pt]{amsart}

\usepackage{amssymb, amsfonts, amsthm, amscd}
\usepackage{epsfig}
\usepackage{graphicx}
\usepackage[all]{xy}
\numberwithin{equation}{section}

\advance\headheight 2pt
\calclayout
\allowdisplaybreaks[3]

\theoremstyle{plain}

\newtheorem{thm}[subsection]{Theorem}
\newtheorem{coro}[subsection]{Corollary}

\newtheorem{lemm}[subsection]{Lemma}

\newtheorem{defn}[subsection]{Definition}

\theoremstyle{definition}

\newtheorem{conj}[subsection]{Conjecture}
\newtheorem{hyp}[subsection]{Hypothesis}

\theoremstyle{remark}

\numberwithin{equation}{section}

\newcommand{\tld}{\tilde}

\newcommand{\s}{\sigma}

\newcommand{\C}{\Bbb C}
\newcommand{\Q}{\Bbb Q}
\newcommand{\Z}{\Bbb Z}
\newcommand{\ra}{\rangle}
\newcommand{\la}{\langle}

\newcommand{\Hom}{\operatorname{Hom}}
\newcommand{\PSL}{\operatorname{PSL}}
\newcommand{\SL}{\operatorname{SL}}
\newcommand{\PGL}{\operatorname{PGL}}
\newcommand{\Ker}{\operatorname{Ker}}

\newcommand{\Syl}{\operatorname{Syl}}

\def\C{{\mathbb C}}

\def\Q{{\mathbb Q}}

\def\Z{{\mathbb Z}}
\def\C{{\mathbb C}}

\makeatother
\makeatletter

\author{Fedor Bogomolov}
\address{Courant Institute of Mathematical Sciences, NYU. \\
 251 Mercer str. \\
 New York, NY 10012, USA}
\email{bogomolo@cims.nyu.edu}

\author{Jorge Maciel}
\address{Courant Institute of Mathematical Sciences, NYU. \\
 251 Mercer str. \\
 New York, NY 10012, USA}
\email{maciel@cims.nyu.edu}

\author{Tihomir Petrov}
\address{Courant Institute of Mathematical Sciences, NYU. \\
 251 Mercer str. \\
 New York, NY 10012, USA}
\email{petrovt@cims.nyu.edu}

\title[Unramified Brauer groups of finite simple groups]
{Unramified Brauer groups of finite simple groups of Lie type $A_{\ell}$}

\begin{document} 



\date{\today}



\begin{abstract}
  We study the subgroup $B_0(G)$ of $H^2(G,{\Q}/{\Z})$ consisting of
  all elements which have trivial restrictions to every Abelian
  subgroup of $G$. The group $B_0(G)$ serves as the simplest
  nontrivial obstruction to  stable rationality of
  algebraic varieties $V/G$ where $V$ is a faithful complex linear
  representation of the group $G$. We prove that $B_0(G)$ is trivial
  for finite simple groups of Lie type $A_{\ell}$.
\end{abstract}

\maketitle

\tableofcontents

\setcounter{section}{1}
\section*{Introduction}
\label{sect:intro}

In this article we answer a group-theoretical question related to the
problem of stable rationality of quotient spaces $V/G$ where $G$ is a
finite (algebraic) group and $V$ is a faithful complex linear
representation of $G$. In algebraic terms, rationality of $V/G$ means
that the field of invariants $\C(V)^G$ is a pure transcendental
extension of the constant field $\C$, and stable rationality means
that there exist a finite number of independent variables
$y_1,\dots,y_k$ such that $\C(V)^G(y_1,\dots,y_k)$ becomes a pure
transcendental extension of $\C$. Simple argument (based on the
geometric version of Hilbert-90 theorem) shows that stable rationality
of $V/G$ does not depend on an individual faithful representation $V$
of $G$ but only on the group $G$ itself. Explicit computations of
invariants show that the varieties $V/G$ are rational in some
nontrivial cases like, for example, the standard representation of the
symmetric group $S_n$.

However, this is not true in general. The first examples of
nonrational and even nonstably rational varieties $V/G$ were obtained
by D. Saltman \cite{Sal}. These solutions of the so-called Noether
problem were obtained by showing that some birational 
invariant, which we denote by $B_0(G)$, is nontrivial
for some series of groups $G$. For any finite group $G$, 
$B_0(G)$ is the subgroup of $H^2(G,{\Q}/{\Z})$ consisting of 
all elements having
trivial restriction on every Abelian subgroup of $G$. It was shown in
\cite{Bog1, Bog2} that $B_0(G)$ coincides with geometric  birational
invariant of a smooth projective model $\widetilde{V/G}$ for $V/G$,
the so-called {\em unramified Brauer group}, introduced earlier by
Artin and Mumford \cite{AM}. Namely,
\[
B_0(G) = {\bf Br}_{\rm nr}(V/G)=H^3(\widetilde{V/G},\Z)_{\rm tors}.
\] 
This fact reduces the computation of the Artin-Mumford invariant 
$V/G$ to a purely group-theoretical question.
 
The group $B_0(G)$ is, in fact, the first of the series of birational
invariants $H^i_{\rm nr}(G)$ of the varieties $V/G$ constructed via
group cohomology (see \cite{CO} for general definitions and
properties, and also \cite{Bog1,Bog3}). In this notation,
$B_0(G)=H^2_{\rm nr}(G)$.
      
However, some results and conjectures in algebraic geometry indicate
that for finite simple groups this kind of birational invariants must
vanish. Our leading hypotheses are the following:
\begin{hyp}
  For any finite simple group $G$ the quotient $V/G$ is stably
  rational.
\end{hyp} 
There is no much evidence to support this statement and, at the
moment, there are no methods to approach this problem except for the
group $A_5=\PGL(2,F_4)$ where it holds.
\begin{hyp} 
  For any finite simple group $G$ the nonramified cohomology groups
  $H^i_{\rm nr}(G)=0$.
\end{hyp}
In this article we test this general hypothesis, formulated in
\cite{Bog3}, for the first nontrivial invariant 
$B_0(G)=H^2_{\rm nr}(G)$ in the case of finite simple groups $G$ of Lie type
$A_{\ell}$ and show that indeed $B_0(G)=0$ for all such groups. Every
group of Lie type $A_{\ell}$ is isomorphic to the 
projective special linear group
$\PSL(n,F_{q})$ over a finite field $F_q$, and the corresponding
second cohomology groups can be found in \cite{Gor}.
 
{\bf Acknowledgments.} We would like to thank Arnaud Beauville for
pointing out an error in the proof of Lemma~\ref{lem41}, part
(2), in the initial version. The first and the third authors were
partially supported by NSF grant DMS-0100837. The second author was
supported by the Instituto Politecnico Nacional - Mexico.

\section{The group $B_{0}(G)$}  
\label{sect:var}

In this section we briefly recall some basic properties of $B_0(G)$.
Let $G$ be a group and $H^{2}(G,\Q/\Z)$ be the set of equivalence
classes of central extensions
\begin{align*}
0\longrightarrow
\mathbb{Q}/\mathbb{Z}\stackrel{i}{\longrightarrow}\tilde{G}\stackrel{j}{\longrightarrow}
G \longrightarrow 1
\end{align*}
of $G$ by $\Q/\Z$. For any $G$-module $\Q/\Z$, elements of
$H^{2}(G,\Q/\Z)$ classify central extensions of the above type. Let
$\mathcal{A}$ be the set of all Abelian subgroups of $G$.
\begin{defn}
The subgroup $B_{0}(G)$ of $H^{2}(G,\Q/\Z)$ is defined by
\begin{align}
  B_{0}(G):= \{\gamma\in H^{2}(G,\Q/\Z)\thickspace\arrowvert
  \thickspace \gamma\mid_{H}=0\thickspace \text{\rm for all}
  \thickspace H\in\mathcal{A}\}.
\end{align}
\end{defn}

\begin{lemm}\label{lm32}
Let $G$ be a group and 
\[
0\longrightarrow\mathbb{Q}/\mathbb{Z}\stackrel{i}{\longrightarrow}\tilde{G}\stackrel{j}{\longrightarrow}
G\longrightarrow 1
\]
be an extension of $G$ by $\Q/\Z$. Let $A$ be an Abelian subgroup of
$G$. For the restriction
\[
0\longrightarrow\mathbb{Q}/\mathbb{Z}\stackrel{i}{\longrightarrow}\tilde{A}\stackrel{j\mid_{\tilde{A}}}{\longrightarrow} A\longrightarrow 1
\]
to be trivial, it is necessary and sufficient that $\tilde{A}$ be an Abelian group.
\end{lemm}
\begin{proof}
  An element $\gamma' \in H^2(A,\Q/\Z)$ is trivial 
  if and only if the corresponding extension 
  $\tilde A$ of $A$ is an Abelian group. It is clear 
  that the trivial extension
  is an Abelian group since it is a direct product of Abelian groups. To
  prove the triviality of the extension $\tilde A$ it is sufficient to
  find a section $s:  A\to \tilde A$. If $x\in A$ has order $n$, then
  there is an element ${\tilde x}\in \tilde A,$ with $j(\tilde x)= x$, of
  order $n$. To find such an element first we take any $y$ with $j(y)=x$.
  Then $y^n\in {\Q/\Z}$. We have the following possibilities
\begin{itemize}
\item If $y^n =0$, then we set $\tilde x:=y$; 
\item If $y^n\ne 0$, then $y^n = a^n$, where $a\in \Q/\Z$, since the group 
  $\Q/\Z$ is infinitely
  divisible. Now we can take $\tilde x := ya^{-1}$. 
\end{itemize}  
Hence for any cyclic subgroup in $A$ we can find a section $s$. Since
$A$ is a direct sum of cyclic groups, the sum of these sections
provides a section for $A$. Therefore, $\tilde A$ is a semi-direct
product of $A$ and $\Q/\Z$, but since it is also a central extension, it
is a direct product. This implies that $\tilde A$ is a trivial central
extension of $A$ if and only if $\tilde A$ is an Abelian group.
Therefore, the fact that the preimage $\tilde{A}=j^{-1}(A)$ is an
Abelian group is equivalent to the triviality of the restriction of
$\gamma \in H^2(G,\Q/\Z)$ on $A$.
\end{proof}

\begin{coro}
  Let $G$ be a group and $\mathcal{X}$ be the set of subgroups of $G$
  generated by two elements whose commutator is the identity element
  of $G$. Then
\[
B_{0}(G)=\{\gamma\in
H^{2}(G,\mathbb{Q}/\mathbb{Z})\thickspace\arrowvert\thickspace\gamma\mid_{H}=0
\thickspace\text{\rm for all}\thickspace H\in\mathcal{X}\}.
\]
\end{coro}

The description above also yields a simple criterion for an element
$\gamma$ of $H^2(G,\Q/\Z)$ not to lie in $B_0(G)$. Let $\tilde G$ be a
central $\Q/\Z$-extension of $G$ defined by $\gamma\in H^2(G,\Q/\Z)$.
Denote by $K_{\gamma}\subset \Q/\Z$ the subgroup in $\Q/\Z$ which lies
in the kernel of every character $\tilde G\to \Q/\Z$. The group
$K_{\gamma}$ is always a finite cyclic group.
\begin{coro}\label{cor35}  
  An element $\gamma$ of $H^2(G,\Q/\Z)$ does not belong to $B_0(G)$ if
  and only if some element of $K_{\gamma}$ can be represented as a
  commutator of a pair of elements $a,b\in \tilde G$.
\end{coro}
\begin{proof}
If $\gamma\in B_0(G)$ then the preimages of any commuting pair of 
elements in $G$ commute in $\tilde G$ as well. Thus if a 
nonzero element $h\in K_{\gamma}$ can be represented as 
$h=ab a^{-1} b^{-1}$, then
the elements $j(a)$ and $j(b)$ commute, but $a$ and $b$ do not. 
This implies that $\gamma$ is not
in $B_0(G)$. Notice that any character of $\tilde G$ is trivial on 
$ab a^{-1} b^{-1}$, and hence any element of the form 
$h=ab a^{-1} b^{-1}$ always belongs to
$K_{\gamma}$.
\end{proof}
\begin{coro}\label{cor35A}
  If the generator of $K_{\gamma}$ is represented as a commutator $ab
  a^{-1} b^{-1}$, then any quotient of $\tilde G$ by a cyclic subgroup
  of ${\Q}/{\Z}$ represents a central
  extension of $G$ which is not in $B_0(G)$.
\end{coro}

\begin{lemm}\label{cor1}
  Let $G$ be a finite group and
\[
H^{2}(G,\mathbb{Q}/\mathbb{Z})=\bigoplus_{p}H^{2}(G,\mathbb{Q}/\mathbb{Z})_{(p)}
\]
be the primary decomposition of $H^{2}(G,\Q/\Z)$, where by
$H^{2}(G,\Q/\Z)_{(p)}$ we denote the $p$-primary component of
$H^{2}(G,\Q/\Z)$.
\begin{enumerate}
\item We have
\[
B_{0}(G)=\bigoplus_{p} B_{0,p}(G),
\]
where $B_{0,p}(G):= B_{0}(G)\cap H^{2}(G,\Q/\Z)_{(p)}$.
\item For every Sylow $p$-subgroup $\Syl_p(G)$ of $G$
  we have an embedding
\[
B_{0,p}(G)\subset B_{0}(\Syl_p(G)).
\]
\end{enumerate}
\end{lemm}
\begin{proof}
  The assertion $(1)$ is a consequence of the fact that a
  central extension of a finite Abelian group is nilpotent and
  therefore decomposes into a direct product of $p$-groups.
  
  The statement in $(2)$ follows from the fact that we have
  an embedding from $H^{2}(G,\Q/\Z)_{(p)}$ into
  $H^{2}(\Syl_p(G),\Q/\Z)$.
\end{proof}
\begin{coro}
  Let $G$ be a finite group and let $\Syl_p(G)$ be a Sylow
  $p$-subgroup of $G$. If $\Syl_p(G)$ is an Abelian group, then
\[
B_{0,p}(G)=0.
\]
\end{coro}
\begin{proof}
  We have $B_{0,p}(G)\subset B_{0}(\Syl_p(G))$ but the latter is zero
  since $\Syl_p(G)$ is Abelian.
\end{proof}
We will also use several results concerning $p$-groups from
\cite{Bog1}.

Let $G$ be a central extension of an Abelian group $A$ with a center
$C=[G,G]$. There is a natural map $\Lambda^2 A\to C$ where 
\[
\Lambda^2 A :=\{x\wedge y \mid x,y\in A\}.
\]
The wedge operation is bilinear, and $x \wedge x=0$.  
It maps $x\wedge y$ to ${\tld x}{\tld y}{\tld x}^{-1}{\tld y}^{-1}$
where ${\tld x}$ and ${\tld y}$ are any two preimages of $x$ and $y$,
respectively, in $G$.
We denote by
$S$ the kernel of the map $\Lambda^2 A\to C$. Obviously, this is a
subgroup of $\Lambda^2 A$. Consider the elements of the form $x\wedge
y \in S$ and denote by $S_{\Lambda}\subset S$ the subgroup in $S$
generated by these elements.
\begin{lemm}
  Let $f:  G\to A$ be a central extension of an Abelian group $A$ by
  $C$. Then
\begin{enumerate}
\item The group $B_0(G)$ is contained in the image
\[  
f^*H^2(A,{\Q}/{\Z})\subset H^2(G,{\Q}/{\Z}).
\]
\item The group $B_0(G)$ is naturally dual to
  $S/S_{\Lambda}$ where $S$ is the kernel of the homomorphism $\Lambda^2
  A\to C$ given by commutators. The group $S_{\Lambda}$ is generated
  by elements of the form $xyx^{-1}y^{-1} \in S$.
\end{enumerate}
\end{lemm}
\begin{proof} 
  See \cite{Bog1}, Lemma 5.1.
\end{proof} 
There are also several results for the more general case of
meta-Abelian groups.
\begin{lemm}\label{lm:38}
  Let $G$ be an extension of a cyclic group by an Abelian group. Then
  $B_0(G)=0$.
\end{lemm}
\begin{proof} 
  See \cite{Bog1}, Lemma 4.9.
\end{proof} 

This result serves as a tool to compute $B_0(G)$ for some meta-Abelian
groups.
 
Let $G$ be a finite meta-Abelian group, $A = G/[G,G]$, and $f: G\to A$
be the canonical projection. Denote by $G_c$ the quotient
$G/[G,[G,G]]$. Let $S$ be the kernel of the linear map $\Lambda^2 A
\to G_c$ and $S_{\Lambda}'$ be the subspace in $S$ generated by
elements $x\wedge y$, where $x,y\in A$, such that there is a pair
$x',y'\in G$ with $f(x')=x, f(y')=y$, and $[x',y']=1$.

The exact sequence of groups
\[
1\longrightarrow [G,G]\longrightarrow G\stackrel{f}{\longrightarrow} A\longrightarrow 1 
\]
defines the spectral sequence
\[
E_2^{p,q}(G) := H^p(A,H^q([G,G],{\Q}/{\Z}))
\]  
which computes the group $H^2(G,{\Q}/{\Z})$. This spectral sequence
provides $H^2(G,{\Q}/{\Z})$ with a filtration
\[
f^*H^2(A,{\Q}/{\Z})\subset K\subset H^2(G,{\Q}/{\Z}), 
\]
so that the quotient $H^2(G,{\Q}/{\Z})/ K$ restricts nontrivially on
$[G,G]$ and $K/f^*H^2(A,{\Q}/{\Z})$ admits an injective map
\begin{align}\label{ar}
r:  K/f^*H^2(A,{\Q}/{\Z})\longrightarrow H^1(A,\Hom([G,G], {\Q}/{\Z})).
\end{align}
This filtration is functorial with respect to the group homomorphisms.
Namely, if $A'\subset A$ is a subgroup then denote the subgroup 
\[
i:  f^{-1}A'\hookrightarrow G, 
\]
as $G'$. Then we have an embedding $G'\hookrightarrow G$, and a natural
cohomology map
\[
i^*:  H^2(G,{\Q}/{\Z})\longrightarrow H^2(G',{\Q}/{\Z}).
\]
This map is a map of filtered groups. Namely, if
\[
1\longrightarrow [G,G]\longrightarrow G'\stackrel{f_1}{\longrightarrow} A'\longrightarrow 1 
\]
is the exact sequence of groups for $G'$, induced from the
corresponding sequence for $G$, then we obtain a natural homomorphism
of spectral sequences
\[
E_2^{p,q}(G)\longrightarrow E_2^{p,q}(G'),
\]  
which is the restriction map 
\[
H^p(A,H^q([G,G],{\Q}/{\Z}))\longrightarrow H^p(A',H^q([G,G],{\Q}/{\Z})).
\]  
This provides the map between cohomology with an invariant 
filtration. Namely, if
\[
f_1^*H^2(A',{\Q}/{\Z})\subset K'\subset H^2(G',{\Q}/{\Z})
\]
is the induced filtration on $H^2(G',{\Q}/{\Z})$, then $i^*$ maps
$K'$ to $K$, and $f^*H^2(A,{\Q}/{\Z})\to f_1^*H^2(A',{\Q}/{\Z})$.

Consider the case of the cyclic subgroup $A' = \Z_m^{(i)}\subset A$.
Denote by $\sigma_i$ the natural restriction map 
\[
\sigma_i:  H^1(A, {\Q}/{\Z})\longrightarrow H^1(\Z_m^{(i)}, \Hom([G,G], {\Q}/{\Z})),
\] 
and by $\sigma$ the sum of such maps over all cyclic subgroups in $A$.
The next result is a weaker version of Theorem 4.2 in \cite{Bog1}
which is better adapted to our setting.

\begin{lemm}\label{lem36}
  Consider the map
\begin{align}\label{sigma}
\sigma:  H^1(A, \Hom([G,G],{\Q}/{\Z}))\longrightarrow \sum H^1(\Z_m , \Hom
  ([G,G],{\Q}/{\Z}))
\end{align}
where $\Z_m\subset [G,G]$ runs through all cyclic subgroups of
$[G,G]$. If $\sigma$ is an embedding, then the group 
$B_0(G)\subset f^* H^2(G,{\Q}/{\Z})$ and 
$B_0(G)=(S/S_{\Lambda}')^*$.
\end{lemm}

\begin{proof} 
  The image $f^* H^2(A,{\Q}/{\Z})$ is identified with the space of
  linear functionals on $S$. By definition, an element $f^* (a)$ belongs
  to $B_0(G)$ if for every Abelian subgroup $B\subset G$ with two
  generators the restriction $f^*(a)\mid_B= 0$. We have to check this
  property only for subgroups $B$ which project onto a subgroup
  with two generators in $A$. Thus the intersection of $B_0(G)$ and
  $f^*H^2(A,{\Q}/{\Z})$ is equal to $(S/S_{\Lambda}')^*$.
  
  The group $B_0(G)$ is contained in $K$ since, by definition, every
  element of $B_0(G)$ restricts trivially on the Abelian group
  $[G,G]$. We want to show that $B_0(G)$ belongs to the kernel of the
  map $r$ on $K$. Assume the contrary, that $\gamma\in K$ and
  $r(\gamma)\neq 0$. If $\gamma\notin \Ker (r)$ then, by assumption on
  the injectivity of $\sigma$, there is a cyclic subgroup
  $\Z^{(i)}_m\subset A$ such that $({\sigma_i\circ r})(\gamma)\ne 0$.
  This implies that the restriction of $\gamma$ on 
  $G'=f^{-1}(\Z_m^{(i)}) \subset G$ is nontrivial, but by Lemma~\ref{lm:38}, the
  group $B_0(G') = 0$.
  
  Therefore, if $\gamma\in B_0(G)$ then $i^*(\gamma)=0$ for any cyclic
  subgroup in $A$. Hence, by our assumption on the $A$-module
  $\Hom([G,G],{\Q}/{\Z})$, we obtain that $B_0(G)\subset \Ker (r)$ and
  therefore it belongs to $f^*H^2(A,\Q/\Z)$.
  
  The image $f^* H^2(A,{\Q}/{\Z})$ is identified with the space of
  linear functionals on $S$. By definition, an element $f^* (a)$
  belongs to $B_0(G)$ if for every Abelian subgroup $B\subset G$ with
  two generators the restriction $f^*(a)\mid_B= 0$. We have to check
  this property only for subgroups $B$ which project onto a subgroup
  with two generators in $A$. Thus the intersection of $B_0(G)$ and
  $f^*H^2(A,{\Q}/{\Z})$ is equal to $(S/S_{\Lambda}')^*$.
\end{proof}   

From Lemma 5.6 in \cite{Bog1} we also easily obtain
\begin{coro}\label{corp6}
  Any $p$-group $G$ with $B_0(G)\ne 0$ has order $\ge p^6$.
\end{coro}

\section{Simple groups and $H^2(G,{\Q}/{\Z})$}

The invariant $B_0(G)\subset H^2(G,{\Q}/{\Z})$ has an analogue,
denoted by $H^i_{\rm nr} (G)$, for group cohomology in any dimension.
For more details we refer to \cite{CO},\cite{CT1},
\cite{Bog1}, and \cite{Bog3}. 
In this notation, $B_0(G)=H^2_{\rm nr}(G)$. 
According to the Bloch-Kato conjecture \cite{BK} all the
cohomology of Galois groups of algebraic closures of functional fields
are, roughly speaking, induced from the Abelian quotients of these
Galois groups. A more geometric version of these conjectures (for
example, see \cite{Bog1, Bog3}) is that, in fact, most of the finite
birationally invariant classes in the cohomology of algebraic
varieties are induced from the similar birational classes of special
$p$-groups. Somehow it looks like finite simple groups do not produce
nontrivial birational invariants of algebraic varieties which leads to
the general hypothesis that all nonramified cohomology of a finite
simple group are trivial (see \cite{Bog1}).

Unfortunately, the computation of general nonramified cohomology
groups is a highly nontrivial task. However, there is a description of
all groups $H^2(G,\Q/\Z)$ for all finite simple groups $G$ (see \cite{Gor},
\cite{Gr}). This fact, together with the rather well
understood structure of Sylow subgroups of finite simple groups, leads
to the following conjecture, which is a particular case of the general
hypotheses stated above:
\begin{conj}[Bogomolov \cite{Bog3}]
  If $G$ is a finite simple group, then
\begin{align}
  H^2_{\rm nr}(G) = B_{0}(G)=0.
\end{align}
\end{conj}

By now the complete list of finite simple groups $G$ is well known. It
is contained, for example, in \cite{Gor}. Apart from a finite number of
sporadic groups all the other are of Lie type. They are parameterized
by a finite number of infinite series depending on the rank of the
corresponding Lie algebra and on the finite base field $F_q$.

Thus any of the simple groups of Lie type has a naturally attached
prime number. The advantage of using a simple group $G$ is that there
is a uniquely defined covering group $\tilde G$. The group $\tilde G$
is a finite universal nontrivial central extension of $G$
with the group $H^2(G,\Q/\Z)$ as a center, so that any central
extension $G'$ of $G$ with $[G',G']=G'$ is a quotient of $\tilde G$.

All the groups $H^2(G,{\Q}/{\Z})$ are also listed (see \cite{Gor}).
The general feature of the corresponding list is that most of the
groups in $H^2(G,{\Q}/{\Z})$ which appear for a given group $G$ in the
series are related to a central extension of the corresponding
algebraic Lie groups. The relevant central extensions have a very
simple description. It is easy to show that they define trivial
elements in $B_0(LG)$ for a Lie group $LG$ over complex numbers (see
\cite{Bog1}, \cite{Bog2}).  Similarly, in each Lie series of finite
groups, apart from a finite number of cases, we have a description of
the corresponding extension in terms of Lie groups.

In this article we prove this result for finite simple groups of Lie
type $A_{\ell}$ where the groups are $\PSL(n,F_q)$. The group
$\SL(n,F_q)$ is a covering group for all $\PSL(n,F_q)$ except for the
following cases: 
\begin{align}\label{pairs}
  (n,q)=(2,4),(2,9),(3,2),(3,4),(4,2).
\end{align}

Therefore, our strategy will be the following:
\begin{itemize}
\item First, we consider the general case when $\SL(n,F_q)$ is the
  universal cyclic extension of $\PSL(n,F_q)$ by the cyclic group
  $\Z_{\gcd(n,q-1)}$ where the group $\Z_{\gcd(n,q-1)}$ is represented
  by diagonal matrices. To prove the result in this case it is
  sufficient to represent the generator $\mu_p$ of any Sylow subgroup
  $\Z_{p^n}$, where $p$ is a prime, of diagonal $\Z_{\gcd(n,q-1)}$
  matrices as a commutator $ABA^{-1}B^{-1}$, where $A,B\in
  \SL(n,F_q)$.
  
\item Second, we consider the above five exceptional pairs
  \eqref{pairs}. Here we do not know the exact description of the
  corresponding central extension. We prove the theorem by
  establishing a stronger result that $B_0(\Syl_p(G)) = 0$ for all
  prime $p$ dividing the order of $H^2(G,{\Q}/{\Z})$.
\end{itemize}

\section{The group ${\SL}(n, F)$ as a central extension}

Let $F$ be a field and let $\mu$ be a primitive $p^n$-th root of unity
in $F$.

\begin{lemm}\label{lem41}
  Assume that $m$ is divisible by $p^n$ where $p$ is prime. The scalar
  matrix $\mu I_{m}$, which belongs to the center of ${\SL}(m, F)$,
  can be written in a commutator form $[A_1,B_1]=A_1B_1A^{-1}_1
  B^{-1}_1$, where $A_1, B_1\in{\SL}(m, F)$.
\end{lemm}

\begin{proof}
We have the following possibilities:
\begin{enumerate}
\item Let $p$ be an odd prime and $A,B$ be two square matrices of
  order $p^n$ over the field $F$ such that $A=(a_{i,j})= (\delta_{i,j}
  \mu^i)$, and $B=(b_{i,j})$ with $b_{i,i+1}=b_{p^n,1}=1$, and
  $b_{i,j}=0$ otherwise. Then $A$ and $B$ belong to
  ${\SL}({p^n},F)$, and $[A,B]=ABA^{-1}B^{-1}=\mu I_{p^n}$.
  
\item If $p=2$ and $n=1$, then the matrices
  $A=\left(\begin{smallmatrix}a&b\\b&-a\end{smallmatrix}\right)$,
  where $a^2 + b^2 = -1$, and $B=\left(\begin{smallmatrix} 0 & 1 \\-1
      & 0\end{smallmatrix}\right)$ satisfy the required condition.
  For a finite field $F_q$, with $q$ odd, such $a$ and $b$ always exist.
  If $F=F_q$ and $q=4k +1$ then $-1$ is a square and we can take
  $a=i$, where $i^2 = -1$, and $b=0$. If $q=4k+3$ then the set of
  all elements of the form $x^2 +y^2$ coincides with $F_q$. Indeed, it
  contains all quadratic residues $x^2$, and it is invariant under the
  multiplication by an arbitrary element $z^2$, where $z\in F_q$.
  Therefore, if it contains at least one non-quadratic element, then
  it coincides with $F_q$. If not, then quadratic residues form an
  additive subgroup in $F_q$. However, $F_q$ does not have an additive
  subgroup of index $2$ for odd $q$. Therefore, $F_q$ is the set of
  elements of the form $x^2 + y^2$. In particular, there are $a$ and
  $b$ such that $a^2 +b^2 = -1$, and they provide entries for $A$.
      
\item If $p=2$ and $n>1$, we select $A=\s_1 X$ and $B =Y \s_2$ where
  $X$ and $Y$ are diagonal matrices and $\s_1,\s_2$ are commuting
  permutations. There are many such pairs. For example, we can
  decompose $2^m$ coordinates into two groups of order $2^{m-1}$ and
  take the diagonal matrix $X$ to be $[1, 1, \dots,i], [1, 1, \dots,
  -i]$ on the main diagonal, where brackets show the boundaries of
  each block. Similarly, take $Y$ to be $[1, \mu, \dots,
  {\mu}^{2^{m-1}-1}], [\mu,{\mu}^2,\dots,{\mu}^{2^{m-1}}]$. The permutation
  $\s_1$ has order $2$ and permutes these two blocks of variables, and
  $\s_2$ permutes variables cyclically within each block. Then $A,
  B\in {\SL}_{2^m}(F)$ and $[A,B]=ABA^{-1}B^{-1}=\mu I_{2^m}$
  
\item If $m$ is divisible by $p^n$, then matrices $A_1$ and $B_1$,
  consisting of $m/p^n$ diagonal blocks of matrices $A$ and $B$
  respectively, also satisfy the relation
  $[A_1,B_1]=A_1B_1A_1^{-1}B_1^{-1} = \mu I_{m }$.
\end{enumerate}
The lemma follows.
\end{proof} 

\begin{coro}
  Any element of the group $H^2(\PSL(n,F), {\Q}/{\Z})$ that defines a
  central extension of $\PSL(n,F)$ which is isomorphic to the quotient
  $\SL(n,F)/{\Z_h}$, where ${\Z_h}$ is a central subgroup, does not
  belong to $B_0(\PSL(n,F))$ (see Corollary~\ref{cor35}). In
  particular, $B_0(\PSL(n,F))=0$ unless the pair $(n,q)$ is
  one of the five exceptional cases in \eqref{pairs}.
\end{coro}

This finishes the proof of our theorem in the general case.

\section{Special simple groups of Lie type $A_{\ell}$}

Consider the remaining cases. These are the following groups:
\begin{enumerate}
\item $G_1=\PSL(2,F_4), H^2(G,{\Q}/{\Z})=\Z_2$
\item $G_2=\PSL(2,F_9), H^2(G,{\Q}/{\Z})=\Z_2\oplus\Z_3$
\item $G_3=\PSL(3,F_2), H^2(G,{\Q}/{\Z})=\Z_2$
\item $G_4=\PSL(3,F_4), H^2(G,{\Q}/{\Z})=\Z_3\oplus \Z_4\oplus \Z_4$
\item $G_5=\PSL(4,F_2), H^2(G,{\Q}/{\Z})=\Z_2$
\end{enumerate}
\begin{lemm}
  In the cases $(1)-(3)$ above the group $B_0(G)=0$.
\end{lemm}
\begin{proof} 
  In the cases $(1)-(3)$ the relevant Sylow $p$-groups have order $<
  p^6$, and therefore, by Corollary~\ref{corp6}, $B_0(G)=0$.
\end{proof}
\begin{lemm} 
  The group $B_0(G_4)=0$.
\end{lemm}
\begin{proof}
  Consider first the element of order $3$ in $H^2(G_4,\Q/\Z)$. The
  corresponding central extension of $G_4= \PSL(3,F_4)$ is the group
  $\SL(3, F_4)$, and hence, by Lemma~\ref{lem41}, the element of order
  $3$ does not belong to $B_0(G_4)$. Similar result for the elements
  order $2$ and $4$ in $H^2(G_4,\Q/\Z)$ follows from the following
  more general result.
\end{proof}

\begin{lemm}
  The group $B_0(\Syl_2(G_4))=0$.
\end{lemm}
\begin{proof} 
  The group $\Syl_2(G_4)$ is a central extension of the Abelian group
  $F_4^2=\Z_2^4 $ by $F_4=\Z_2^2$. Thus we again can use the general
  formula $B_0(\Syl_2(G_4))=(S/S_{\Lambda})^*$. We have two generators
  $\la x,y\ra$ over $F_4$, which is linearly generated over $\Z_2$ by
  $\la 1, s\ra$ with $s^3 = 1$. The commutator map $\Lambda^2 \Z_2^4
  \to F_4$ coincides with the map between the second exterior power of
  $\Z_2^4$, considered as $\Z_2$-space, and the second exterior power
  of the same space, considered as $F_4$-space. Thus $S_{\Lambda}$ is
  generated by elements $[u,su]$ where $u\in \Z_2^4=F_4^2$. The space
  $S$ has dimension $4$ over $\Z_2$ and a basis $\la x,y, sx, sy\ra$.
  To finish the proof we will show that the elements
\[
[x,sx],\ [y,sy],\ [x+y,s(x +y)],\ [x+sy,s(x +sy)]
\]
are linearly independent over $\Z_2$. Indeed, this set is linearly
equivalent to 
\[
[x,sx],\ [y,sy],\ ([x,sy] + [y ,sx]),\ ([x,y]+[x,sy]+[sy,sx]), 
\]
and the last one is linearly independent over $\Z_2$. Therefore,
$S_{\Lambda}=S$ and $B_0(\Syl_2(G_4))=0$, which implies that
$B_0(G_4)=0$.
\end{proof}    

\begin{thm}\label{tm63}
  The group $B_0(\Syl_2(G_5))=0$.
\end{thm}
\begin{proof} 
  The group $\Syl_2(G_5)$ is a meta-Abelian $2$-group. We are going to
  use Lemma~\ref{lem36} in order to compute the group
  $B_0(\Syl_2(G_5))$.
  
  The group $\Syl_2(G_5)$ is the group of upper-triangular
  $\Z_2$-matrices. It is generated by the elements $\la a_{i,i+1}\ra$,
  with $i=1,2,3$. Here $(a_{i,j})$ denotes the matrix having $1$'s
  on the main diagonal and $0$'s everywhere else apart from the
  $(i,j)$'th entry, which is equal to $1$. Denote the group
  $\Syl_2(G_5)$ by $G$ to simplify the notations. The group $[G,G] = \Z_2^3$
  is generated  by the elements $\la a_{1,3}, a_{1,4},
  a_{2,4}\ra$.
\begin{lemm}\label{lm64} 
  The conditions of Lemma~\ref{lem36} are satisfied for the group $G =
  \Syl_2(G_5)$.
\end{lemm}
\begin{proof}
  The module $\Hom([G,G],{\Q}/{\Z})=\Z_2^3$. For simplicity, we will
  denote it by $W$. Denote the generators of the quotient $A=G/[G,G]$
  by $a_{i,i+1}, i=1,2,3$, and the generators of $W$ by $a_{1,3}^*,
  a_{1,4}^*, a_{2,4}^*$, respectively.

The group $A$ acts on $W$. The action of $A$ is trivial on 
$V=\la a_{1,3}^*,a_{1,4}^*\ra\subset W$, 
and on the quotient
$W/V=C$. We have an exact sequence
\[
H^0(A,C)\stackrel{\delta}{\longrightarrow} H^1( A, V)\longrightarrow
H^1( A, W)\longrightarrow H^1 (A,C).
\]
Since the action of $A$ is trivial on $V$ and $C$, any element in
$H^1(A,V)$ and $H^1(A,C)$ restricts nontrivially on some cyclic
subgroup of $A$.
  
Now we have to prove the same for $H^1(A,V)/\delta H^0(A,C)$. The
group $H^0(A,C)=\Z_2$, and we denote the only nontrivial element in
this group by $w$. Denote the corresponding nontrivial element of
$H^1(A,V)$ by $\delta w$. We have that $A=\Z_2^3, V=\Z_2^2$ and
$H^1(A,V)=A^*\otimes V$. Thus we can select generators 
$\la x_1, x_2, x_3 \ra$ of $A$ and $\la v_1,v_2\ra$ of $V$ so that either
\begin{itemize}
\item $\delta w = x_1^* \otimes v_1$, or
\item $\delta w = x_1^* \otimes v_1 + x_2^* \otimes v_2$.
\end{itemize}
Here $\la x^*_1, x^*_2, x^*_3\ra$ is a basis for $A^*$, such that 
$x_i^*(x_j)=\delta_{i,j}$.
  
In order to finish the proof, we need the following lemma
\begin{lemm}
  Let $y\subset A^* \otimes V$ be a nonzero element such that $y\neq \delta w$.
  Then there is a projection 
\begin{align}\label{proj}
  p_y:  A^*\longrightarrow \Z_2,
\end{align}
such that the images of the elements $\delta w$ and $y$ under the
induced projection $A^*\otimes V \to \Z_2\otimes V$ are different, and
the image of $y$ is nonzero.
\end{lemm}
\begin{proof}
Let $y =\sum b_{i,j} x_i^*\otimes v_j$.
\begin{itemize}
\item If $b_{i,j} \neq 0$ for $(i,j)\neq(1,1)$, then the
  projection of $y$ on one of the groups 
  $\Z_2 x_1^*\otimes V, \Z_2  x_2^*\otimes V,$ or $\Z_2 x_3^*\otimes V$ 
  is nonzero and different from
  $\delta w$. Otherwise, $y = \delta w$ which proves the result in the
  first case.
  
\item If $b_{3,j}\neq 0$ for some $j$, then the restriction of $y$ on
  $\Z_2 x_3^*\otimes V$ is nontrivial and nonequal to $w$.
  
\item If either $b_{1,2}$ or $b_{2,1}$ is nonzero, then the
  restriction of $y$ on either 
  $\Z_2 x_1^*\otimes V$ or $\Z_2 x_2^*\otimes V$ is nonzero and
  different from the corresponding restriction of $\delta w$.
\end{itemize}  
Therefore, we can assume that either 
\begin{itemize}
\item $y=x_1^*\otimes v_1$, or
\item $y=x_2^*\otimes v_2$. 
\end{itemize}
In both cases, the restriction of $y$ on $\Z_2 (x_1^* + x_2^*)$ will be
nonzero and different from the restriction of $\delta w$.
\end{proof}
Thus if we take the subgroup $p_y^* :  \Z_2\hookrightarrow A$, 
where $p_y^*$ is dual to the map in \eqref{proj}, then 
the restriction of $y$ on $p_y^* \Z_2$ is nonzero 
in $H^1(\Z_2,W)$. This finishes the proof of Lemma~\ref{lm64}.
\end{proof}
    
\begin{coro}
  By Lemma~\ref{lem36}, we obtain that $B_0(\Syl_2(G_5))$ is contained
  in the image $f^* H^2(A,\Q/\Z)$ where $A$ is the maximal Abelian
  quotient of $\Syl_2(G_5)$.
\end{coro}

Now, in order to finish the proof of the triviality of $B_0(G)$ we
have to compute the groups $f^*H^2(A,\Q/\Z)$ and $S$ in the above
case. Denote the basis $\Z_2$-characters of the group $\Syl_2(G_5)$
by $\la a_{i,i+1}^*\ra$, $i=1,2,3$.

\begin{lemm}
  The group $f^*H^2(A,\Q/\Z)= \Z_2$ and it is generated by an element
  $a_{1,2}^*\wedge a_{3,4}^*$ which is nontrivial on the commutative
  group generated by $\la a_{1,2}, a_{3,4} \ra$.
\end{lemm}

\begin{proof}
  The group $H^2(A,\Q/\Z)= \Z_2^3$. Let us show that both elements
  $a_{1,2}^*\wedge a_{2,3}^*$ and $a_{2,3}^* \wedge a_{3,4}^*$ belong
  to the kernel of $f^*$. Consider two projections of
  $\Syl_2(G_5)$ onto a central extension of $\Z_2\oplus\Z_2$, which is
  the upper-triangular group of $3\times 3$ matrices. The first one is
  given by deleting the fourth column, and the second one by deleting
  the first row. This shows that the elements $a_{1,2}^*\wedge
  a_{2,3}^*$ and 
  $a_{2,3}^* \wedge a_{3,4}^*$ are trivial already on the
  corresponding triangular groups, and hence they are trivial on
  $\Syl_2(G_5)$. The group $\la
  a_{1,2}, a_{3,4} \ra$ is an Abelian $B=\Z_2\oplus\Z_2$-subgroup of
  $\Syl_2(G_5)$. The element $f^*(a_{1,2}^*\wedge a_{3,4}^*)$ is
  nontrivial on $B$ which shows that $a_{1,2}^*\wedge a_{3,4}^*\neq
  0$, and it does not belong to $B_0(\Syl_2(G_5))$.
\end{proof}   
 
Thus we proved in particular that $B_0(G_5)=0$ which finishes the
proof of Theorem~\ref{tm63}. 
\end{proof}
Simultaneously we finished the proof of our main result.
  
\begin{thm}[Main]
  For every finite field $F_{q}$ of order $q$ and every integer
  $n\geqslant 2$,
\begin{align*}
  B_{0}(\PSL(n,F_{q}))=0.
\end{align*}
\end{thm}

\end{document}